\documentclass{amsart}
\usepackage{amssymb,amsmath,amsthm,graphics,amscd,amsfonts}

\usepackage{color}

\theoremstyle{plain}
\newtheorem{theorem}{Theorem}[section]

\newtheorem{corollary}[theorem]{Corollary}
\newtheorem{lemma}[theorem]{Lemma}

\newtheorem{definition}[theorem]{Definition}

\newtheorem{fact}[theorem]{Fact}

\newcommand{\bfC}{{\mathbb C}}
\newcommand{\bfP}{{\mathbb P}}
\newcommand{\bfR}{{\mathbb R}}

\newcommand{\bfQ}{{\mathbb Q}}

\newcommand{\barj}{{\overline j}}

\newcommand{\barz}{{\overline z}}

\newcommand{\barpartial}{{\overline \partial}}

\newcommand{\mapright}[1]{\smash{\mathop{   \hbox to 0.7cm{\rightarrowfill}}
  \limits^{#1}}}

\def\bp{\overline{\partial}}

\def\om{\omega}

\def\p{\partial}
\def\bp{\overline{\partial}}

\def\wt{\widetilde}

\begin{document}

\title{Multiplier ideal sheaves and geometric problems}
\author{Akito Futaki}
\address{Department of Mathematics, Tokyo Institute of Technology, 2-12-1,
O-okayama, Meguro, Tokyo 152-8551, Japan}
\email{futaki@math.titech.ac.jp}
\author{Yuji Sano}
\address{Department of Mathematics, Kyushu University,
6-10-1, Hakozaki, Higashiku, Fukuoka-city, Fukuoka 812-8581 Japan}
\email{sano@math.kyushu-u.ac.jp}

\date{July 26, 2009}

\begin{abstract} 
In this expository article we first give an overview on multiplier ideal sheaves and
geometric problems in K\"ahlerian and Sasakian geometries.
Then we review our recent results on the relationship between the support of the subschemes cut out by multiplier ideal sheaves 
and the invariant whose non-vanishing obstructs the existence of K\"ahler-Einstein metrics on
Fano manifolds.
\end{abstract}



\maketitle

\section{Introduction}
One of the main problems in K\"ahlerian and Sasakian geometries is the existence problem of Einstein metrics. 
An obvious necessary condition for the existence of a K\"ahler-Einstein metric on a compact K\"ahler manifold $M$ is that the first Chern class
$c_1(M)$ is negative, zero or positive since the Ricci form represents the first Chern class. 
This existence problem in K\"ahlerian geometry was settled by
Aubin \cite{aubin76} and Yau \cite{yau78} in the negative case and by Yau \cite{yau78} in the zero case. In the remaining case when the manifold has 
positive first Chern class, in which case the manifold is called a Fano manifold in algebraic geometry, there are two known
obstructions. One is due to Matsushima \cite{matsushima57} which says that the Lie algebra $\mathfrak h(M)$ of all holomorphic vector fields on a compact
K\"ahler-Einstein manifold $M$ is reductive. and the other one is due to the first author \cite{futaki83.1} which is given by a Lie algebra character $F : \mathfrak h(M) \to \bfC$
with the property that if $M$ admits a K\"ahler-Einstein metric then $F$ vanishes identically.
Besides, it has been conjectured by Yau \cite{yau93} that a more subtle condition related to geometric invariant theory 
should be equivalent to the existence of K\"ahler-Einstein metrics. This idea was made explicit in the paper \cite{tian97} of Tian  in which a notion called K-stability
was introduced. Tian used a generalized version of the invariant
$F$ for normal almost Fano varieties and used it as the numerical invariant for the stability condition. 
 The link between the idea of GIT stability and geometric problems such as the existence problems of Hermitian-Einstein metrics and constant scalar curvature K\"ahler metrics
 can be explained through the moment maps in symplectic geometry. 
 The explanation from this viewpoint can be found for example 
in \cite{donkro}, \cite{fujiki92}, \cite{donaldson97}. 
Recall that an extremal K\"ahler metric is by definition a K\"ahler metric such that the gradient vector field of the scalar curvature is a holomorphic vector field.
In particular, a K\"ahler metric of constant scalar curvature is an extremal K\"ahler metric.
The theorem of Matsushima is extended for extremal K\"ahler manifolds by Calabi \cite{calabi85} 
as a structure theorem of the Lie algebra $\mathfrak h(M)$ 
on an extremal K\"ahler manifold $M$,
and the first author's obstruction $F$ can be extended as an obstruction to the existence of constant scalar 
curvature K\"ahler metric in a fixed K\"ahler class (\cite{futaki83.2}, \cite{calabi85}).
The theorem of Calabi and the character $F$
are explained in the framework of the moment map picture by X. Wang \cite{xwang04} (see also \cite{futaki05}).

In \cite{donaldson02} Donaldson refined the notion of K-stability for 
a polarized manifold $(M, L)$, that is, a pair of an algebraic manifold $M$ and 
an ample line bundle $L$ over $M$,  and  conjectured that
there would exist a K\"ahler form in $c_1(L)$ of constant scalar curvature if and only if $(M,L)$ is K-polystable. To define K-(poly)stability for $(M, L)$ Donaldson refined the invariant
$F$ even for non-normal varieties which are degenerations of the polarized manifold $(M,L)$ and used it as the numerical invariant for the stability condition. 
The K-stability is defined as follows. 
For an ample line bundle $L$ over a projective variety  $M$ of dimension $m$, a test configuration of
exponent $r$ consists of the following. \\
(1)\ \ A flat family of schemes $\pi : {\mathcal M} \to \bfC$:\\
(2)\ \ $\bfC^*$-action on ${\mathcal M}$ covering the usual $\bfC^*$-action on $\bfC$:\\
(3)\ \ $\bfC^*$-equivariant line bundle $\mathcal{L} \to {\mathcal M}$ such that
\begin{itemize}
\item 
for $t \ne 0$ one has $M_t = \pi^{-1}(t) \cong M$ and
$(M_t, {\mathcal L}|_{M_t}) \cong (M, L^r)$,
\item 
$\chi(M_t, L^r_t) = \sum_{p=0}^m (-1)^p \dim 
H^p(M_t, L_t^r)$ does not depend on $t$, in particular for $r$ sufficiently large
$\dim H^0(M_t, L_t^r) = \dim H^0(M,L^r)$ for all $t \in \bfC$. Here we write $L^r_t$ for ${\mathcal L}|_{M_t}$
though $L$ may not exist for $t=0$.
\end{itemize}

\noindent
The $\bfC^*$-action on $(\mathcal{L}, {\mathcal M})$ induces a $\bfC^*$-action on the central fiber
$L_0 \to M_0 = \pi^{-1}(0)$. Moreover if 
$(M,L)$ admits a $\bfC^*$-action, then one obtains a test configuration
by taking the direct product $M \times \bfC$. This is called a product configuration.
A product configuration is called a trivial configuration if the action of $\bfC^*$ on $M$ is trivial.

\begin{definition}\ \ $(M,L)$ is said to be K-semistable (resp. stable) if
the invariant $F_1$ (defined below) of the central fiber $(M_0, L_0)$ is non-positive (resp. negative)
for all non-trivial test configurations. $(M,L)$ is said to be
K-polystable if it is K-semistable and $F_1 = 0$ only
if the test configuration is product. 
\end{definition}
\noindent
Here the invariant $F_1$ is defined as follows.
Let 
$L_0 \to M_0$ be an ample line bundle over an $m$-dimensional projective 
scheme. We assume that a ${\mathbb C}^*$-action as bundle isomorphisms of 
$L_0$ covers the ${\mathbb C}^*$-action on $M_0$. 
For any positive integer $k$, there is an induced $\bfC^*$ action on
$W_k = H^0(M_0, L_0^k)$. Put
$d_k = \dim W_k$ and let $w_k$ be the weight of $\bfC^*$-action on 
$\wedge^{d_k}W_k$. 
For large $k$, $d_k$ and $w_k$ are polynomials in $k$ of degree $m$  and $m+1$ respectively
by the Riemann-Roch and the equivariant Riemann-Roch theorems. For sufficiently large $k$ we expand
$$ \frac{w_k}{kd_k} = F_0 + F_1k^{-1} + F_2k^{-2} + \cdots. $$
When $M_0$ is smooth, $F_1$ coincides with  $F(X)$ up to a negative multiple constant where $X$ is the infinitesimal generator of the $\bfC^{\ast}$-action on $M_0$.

The necessity of K-polystability for the existence of constant scalar curvature K\"ahler metric has been studied by Chen and Tian \cite{chentian04}, Paul and Tian \cite{paultianCM1}, 
\cite{paultianCM2}, Donaldson \cite{donaldson05}, Stoppa \cite{stoppa0803} and
Mabuchi \cite{mabuchi0812}.

Returning to Fano manifolds, there are two hopeful approaches to prove the sufficiency of K-polystability for the existence of K\"ahler-Einstein metrics. One is the Monge-Amp\`ere equation and
the other is the K\"ahler-Ricci flow. In both cases the difficulty arises in the $C^0$-estimate, and when the $C^0$-estimate fails the multiplier ideal sheaves and the subschemes
cut out by them appear. 
The multiplier ideal sheaves arising from the Monge-Amp\`ere equation were studied by Nadel \cite{nadel90}, and those arising from Ricci flow were studied for example  by 
Phong, Sesum and Sturm \cite{pss0611} and Rubinstein \cite{rubinstein0708}. We will give an overview on this subject in section 2.
On the hand, for a given subscheme $V$ in $M$, Ross and Thomas \cite{rossthomas06} considered the test configuration obtained by
blowing up $M \times \bfC$ along $V \times \{0\}$. $M$ is said to be {\it slope stable} if the invariant $F_1$ for the test configuration is negative for any $V$. 
These works lead us to ask how  the invariant $F$ (or more generally $F_1$)  is related to the multiplier ideal sheaves arising from the Monge-Amp\`ere equation
and the Ricci flow. We will treat this subject in section 3.

Now we turn to Sasakian geometry. 
For general facts about Sasakian geometry, refer to the book \cite{boyer-galicki08}.
Let $(S,g)$ be a Riemannian manifold. We denote its Riemannian cone
$(\mathbb{R}_+\times S,dr^2+r^2g)$ by $(C(S),\bar{g})$.
A Riemannian manifold $(S,g)$ is said to be a {\bf Sasaki manifold} if 
the Riemannian cone $(C(S),\overline{g})$ is K\"ahler.
From this definition the dimension of Sasaki manifold $(S,g)$ is
odd, and we put $\dim_{\bfR} = 2m + 1$ so that $\dim_{\bfC} C(S) = m+1$.
$(S,g)$ is isometric to the submanifold $\{r = 1\} = \{1\} \times S 
\subset (C(S),\overline{g})$, and we identify $S$ with the submanifold $\{r = 1\} $. 
Let $J$ be the complex structure on 
$C(S)$ giving the K\"ahler structure.
Consider the vector field 
$$\wt{\xi}=Jr\frac{\partial }{\partial r}.$$
Then $\frac12 (\wt{\xi} - iJ\wt{\xi})$ is a holomorphic vector field. The restriction $\xi$  of $\wt{\xi}$ to $S \cong \{r=1\}$ becomes a
Killing vector field, called the Reeb vector field. The flow generated by $\xi$ is called the Reeb flow.
The restriction $\eta$ of the 
$1$-form $\wt{\eta}$ on $C(S)$ defined as
$$\wt{\eta}=\frac{1}{r^2}\bar{g}(\Tilde{\xi},\cdot)
=\sqrt{-1}(\bar{\partial}-\partial)\log r$$
to $S \cong \{r=1\}$ becomes a contact form. Hence $d\eta$ defines K\"ahler forms on local orbit
spaces of Reeb flow.
That is to say, the $1$-dimensional foliation defined by $\xi$ comes equipped with a structure of transverse K\"ahler foliation.
A Sasaki manifold is said to be regular if the Reeb flow generates a free $S^1$-action, quasi-regular if all the orbits are closed.
A Sasaki manifold is said to be irregular if it is not quasi-regular.

For a polarized manifold $(M,\omega)$ the associated $U(1)$-bundle $S$ of $L$ becomes a regular Sasaki manifold in a natural way:
Choose a positive $(1,1)$ form representing $c_1(L)$, take the Hermitian metric $h$ on $L$ such that 
the connection form $\wt\eta$ on $L$ has its curvature form $d\wt\eta$ equal to $\omega$. The K\"ahler cone $C(S)$ is $L$ minus the zero section 
with the K\"ahler form given by $\frac i2\p\bp r^2$ where $r$ is the distance from the zero section.
Conversely, any regular Sasaki manifold is given in this way. Similarly, a quasi-regular Sasaki manifold is given as an associated $U(1)$-orbibundle over an orbifold. 

As is shown in \cite{FOW} most of ideas in K\"ahler geometry can be extended to transverse K\"ahler geometry for Sasaki manifolds. For example one can extend Calabi's theorem
to compact Sasaki manifold whose transverse K\"ahler metric is an extremal K\"ahler metric, and one can extend the obstruction $F$ as an obstruction for a basic cohomology class
to admit a transverse K\"ahler form with constant scalar curvature. 

A Sasaki-Einstein manifold is a Sasaki manifold whose metric is an Einstein metric. This condition is equivalent to that the transverse K\"ahler metric is K\"ahler-Einstein.
Thus the study of the existence problem of Sasaki-Einstein metrics are closely related to the problem of K\"ahler-Einstein metrics.
But there are differences between them. To explain the differences let $k$ be the maximal dimension of the torus which acts on $C(S)$ as holomorphic isometries. When $k = m+1$
the cone $C(S)$ is a toric variety, and in this case the Sasaki manifold $S$ is said to be toric. Notice that $k$ is at least $1$ because $\wt\xi$ generates holomorphic isometries on $C(S)$.
The other extreme case is therefore when $k = 1$. In this case the Sasaki manifold is necessarily quasi-regular. 

The contact bundle $D = \mathrm{Ker}\,\eta \subset TS$ has a complex structure given by the restriction of $J$. A necessary condition for the existence of Sasaki-Einstein metric
in a fixed transverse K\"ahler structure is that the following two conditions are satisfied:\\
(a)\ the basic first Chern class is represented by a positive transverse $(1,1)$-form;
(b)\  $c_1(D) = 0$, \\
see \cite{boyer-galicki08} or \cite{FOW} for the proof. 
In \cite{FOW} it is proved that if a compact toric Sasaki manifold satisfies the conditions (a) and (b)  then we can deform the Reeb vector field so that the
resulting Sasaki manifold has a Sasaki-Einstein metric. It can be shown that the conditions (a) and (b) are rephrased as the Sasaki manifold is obtained from
the toric diagram of  constant height, and equivalently as the apex of $C(S)$ is a $\bfQ$-Gorenstein singularity (c.f. \cite{CFO}).
In the case of the other extreme when $k = 1$ the conditions (a) and (b) only say that the orbit space of the Reeb flow is a Fano orbifold. In this case there is no deformation space
of Reeb vector field and the problem has the same difficulty as the problem of K\"ahler-Einstein metrics. For the intermediate cases when $1 < k < m+1$ the authors do not even know how to
state the conjecture. In the extreme case where $k = 1$ numerous existence results were obtained by Boyer, Galicki, Koll\'ar and their collaborators using the multiplier ideal
sheaves, which will be reviewed in the next section.

\section{An overview of multiplier ideal sheaves}
In this section, we recall the results about the relationships between the existence of K\"ahler-Einstein metrics on Fano manifolds and the multiplier ideal sheaves, and related topics.
In particular we focus on Nadel's works and recent results about the multiplier ideal sheaves and the K\"ahler-Ricci flow.

Nadel \cite{nadel90} gave a sufficient condition for the existence of K\"ahler-Einstein metrics on Fano manifolds by using the multiplier ideal sheaves, which was originally studied in the works of J.J. Kohn.
Let $M$ be a compact $m$-dimensional Fano manifold.
Let $g$ be a K\"ahler metric on $M$, whose K\"ahler class equals to the first Chern class $c_1(M)$ of $M$.
Let $\gamma_0\in (0,1)$.
We denote the K\"ahler form and the Ricci form of $g$ by $\omega_g$ and $\mbox{Ric}(g)$ respectively.
Let 
\begin{equation}\label{eq:seq_potentials}
	S:=\{\varphi_k\in C^\infty_\mathbb{R}(M) \mid
	g_{i\bar{j}}+\partial_i\partial_{\bar{j}}\varphi_k>0, \,\,
	\sup_M\varphi_k=0, \,\,
	1\le k <\infty
	\}
\end{equation}
be a sequence of K\"ahler potentials with respect to $\omega_g$ such that 
\begin{equation}\label{eq:integral_condition}
	\lim_{k\to\infty}\int_M e^{-\gamma\varphi_k} dV = \infty
\end{equation}
for any $\gamma \in (\gamma_0, 1)$, and that there is a nonempty open subset $U\subset M$ satisfying that 
\begin{equation}\label{eq:condition2}
	\int_U e^{-\varphi_k}dV \le O(1)
\end{equation}
as $k\to \infty$, where $dV$ is a fixed volume form.
Remark that the last condition (\ref{eq:condition2}) always holds for any $S$ due to $g_{i\bar{j}}+\partial_i\partial_{\bar{j}}\varphi>0$ (see for instance \cite{tian87}).
For each $S$, Nadel constructed a coherent ideal sheaf $\mathcal{I}(S)$, which is called the multiplier ideal sheaf (MIS).
We will explain later the simpler definition of multiplier ideal sheaves given by Demailly-Koll\'ar \cite{demailly-kollar01}.

Here let us recall the outline of Nadel's construction. (See the original paper \cite{nadel90} for the full details.)
Let $L$ be an arbitrary ample line bundle on $M$ which is not necessarily the anticanonical line bundle of $M$.
We define $H^0(M,L^\nu)_S$ to be the set of all $f\in H^0(M,L^\nu)$ for which there exists a sequence $\{f_k\}$ of $H^0(M, L^\nu)$ such that
\[
	\int_M |f_k|^2e^{-\gamma\varphi_k}dV \le C
\]
for some $\gamma\in (\gamma_0,1)$ and $f_k\to f$ uniformly.
Consider the homogeneous coordinate ring
\[
	R(M,L)=\bigoplus_{\nu=0}^{\infty}H^0(M,L^\nu).
\]
Define 
\[
	I(M,S,L)=\bigoplus_{\nu=0}^{\infty}H^0(M,L^\nu)_S,
\]
which is a homogeneous ideal $I(M,S,L)$ of the graded ring $R(M,L)$.
Then, the ideal sheaf $\mathcal{I}(S)$ is defined as the algebraic sheaf of ideals on $M$ associated to $I(M,S,L)$.
It is proved in \cite{nadel90} that this construction is independent of the choice of $L$.
Let $\mathcal{V}(S)$ be the (possibly non-reduced) subscheme in $M$ cut out by $\mathcal{I}(S)$. 
This subscheme is characterized as follows.
A point $p\in M$ is contained in the complement of $\mathcal{V}(S)$ if and only if there exist an open neighborhood $W$ of $p$ in $M$ and a real number $\gamma\in (\gamma_0,1)$ such that 
\[
	\int_W e^{-\gamma\varphi_k}dV \le O(1)
\]
as $k\to \infty$.
Remark that $\mathcal{V}(S)$ is neither empty nor $M$ if $S$ satisfies the conditions (\ref{eq:seq_potentials}), (\ref{eq:integral_condition}) and (\ref{eq:condition2}).

One of the distinguished properties of this ideal sheaf is the following vanishing theorem.
\begin{theorem}\cite{nadel90}
For every semi-positive Hermitian holomorphic line bundle $L$ on $M$,
\[
	H^i(M, \mathcal{O}(L)\otimes \mathcal{I}(S))=0,\,\,
	\mbox{for all } i>0. 
\]
\end{theorem}
In particular, we have
\begin{equation}\label{eq:vanishing_MISheaf}
	H^i(M, \mathcal{I}(S))=0, \,\, i \ge 0.
\end{equation}
Remark that (\ref{eq:vanishing_MISheaf}) at $i=0$ follows from the fact  that the subscheme $\mathcal{V}(S)$ is not empty.
We also find that  (\ref{eq:vanishing_MISheaf}) implies 
\begin{equation}	\label{eq:vanishing_MIScheme}
\left.\begin{array}{cccc}
H^0(\mathcal{V}(S),\mathcal{O}_{\mathcal{V}(S)}) & = & \mathbb{C} &
\\
H^i(\mathcal{V}(S),\mathcal{O}_{\mathcal{V}(S)}) & = & 0 &\mbox{for all }i>0.
\end{array}\right.
\end{equation}
This vanishing formula (\ref{eq:vanishing_MIScheme}) gives us several geometric properties of $\mathcal{V}(S)$.
For example, 
\begin{enumerate}
\item
$\mathcal{V}(S)$ is connected.
\item\label{list:zero_dim}
If $\mathcal{V}(S)$ is zero dimensional, then it is a single reduced point.
\item
If $\mathcal{V}(S)$ is one dimensional, then it is a tree of smooth rational curves.
\end{enumerate}
The main result in \cite{nadel90} is that if a Fano manifold does not admit K\"ahler-Einstein metrics then the bubble of the solution of the continuity method induces a proper multiplier ideal sheaf with the above vanishing formula (\ref{eq:vanishing_MIScheme}).
To explain it, let us recall the continuity method for the Monge-Amp\`ere equation.
Here, we assume that $\gamma_0=m/(m+1)$.
Consider the following equation
\begin{equation}\label{eq:Monge-Ampere}
	(\det(g_{i\bar{j}})+\partial_i\partial_{\bar{j}}\varphi_{t})/
	(\det(g_{i\bar{j}}))=
	\exp(h_g-t\varphi_t),
\end{equation}
where $t\in[0,1]$ and $h_g$ is the real-valued function defined by

\begin{equation}\label{eq:ricci_potential}
	\mbox{Ric}(g)-\omega_g=\frac{\sqrt{-1}}{2\pi}\partial\bar{\partial}h_g,
	\,\,
	\int_M e^{h_g} \omega_g^m=\int_M \omega_g^m=V.
\end{equation}
It is well-known that the space $T:=\{t\in [0,1] \mid \mbox{(\ref{eq:Monge-Ampere}) has a solution}\}$ contains $0$ (due to the Calabi-Yau theorem) and open in $[0,1]$ (due to the implicit function theorem).
If $T$ is closed then (\ref{eq:Monge-Ampere}) is solvable at $t=1$, i.e., $\omega_g+\frac{\sqrt{-1}}{2\pi}\partial\bar{\partial}\varphi_{1}$ gives a K\"ahler-Einstein form.
A priori estimates for the closedness of $T$ were given by Yau \cite{yau78};
if (\ref{eq:Monge-Ampere}) is solvable at $s\in [0,t)$ and $\|\varphi_s\|_{C^0} $ is uniformly bounded then  (\ref{eq:Monge-Ampere}) is solvable at $s=t$.
Nadel proved that if the solution $\{\varphi_t\}_{0\le t < t_0}$ of (\ref{eq:Monge-Ampere}) violates the above estimate, then
there is a sequence $\{t_k\}$ such that $t_k\to t_0$ as $k\to\infty$ and $\{\varphi_{t_k}-\sup_M\varphi_{t_k}\}_{k=1}^\infty$ induces a proper multiplier ideal sheaf $\mathcal{I}$.
In this paper, we call it the \textbf{ K\"ahler-Einstein multiplier ideal sheaf} (KE-MIS). 
Summing up,
\begin{theorem}[\cite{nadel90}]\label{thm:Nadel_MIS}
Let $M$ be a Fano manifold which does not admit K\"ahler-Einstein metric.
Let $G$ be a compact subgroup of the group $ \mbox{Aut}(M)$ of holomorphic automorphisms of $M$.
Assume that $M$ does not admit any $G^\mathbb{C}$-invariant proper multiplier ideal sheaf. Then $M$ admits K\"ahler-Einstein metrics. Here $G^\mathbb{C}$ denotes the complexification of $G$.
\end{theorem}
By combining the above theorem and the geometric properties of $\mathcal{V}(S)$ given by the vanishing formula (\ref{eq:vanishing_MIScheme}), Nadel gave many examples of K\"ahler-Einstein Fano manifolds.
Recently Heier \cite{heier0710ii} applied this method to (re-)prove the existence of K\"ahler-Einstein metrics on complex del Pezzo surfaces obtained from the blow up of $\mathbb{CP}^2$ at 3,4 or 5 points, which was originally proved by Siu \cite{siu88}, Tian \cite{tian87}, Tian and Yau \cite{tian-yau87}.

This method was extended to the case of Fano orbifolds by Demailly-Koll\'ar \cite{demailly-kollar01}.
Their construction  is simpler  than \cite{nadel90}.
Let $\psi$ be an $\omega_g$-plurisubharmonic (psh) function (or almost psh function with respect to $\omega_g$), i.e., a real-valued upper semi-continuous function satisfying $\omega_g+\frac{\sqrt{-1}}{2\pi}\partial\bar{\partial}\psi\ge 0$ in the current sense.
The multiplier ideal sheaf with respect to $\psi$ in the sense of \cite{demailly-kollar01} is  the ideal sheaf defined by the following presheaf
\begin{equation}\label{eq:MIS_PSH}
	\Gamma(U, \mathcal{I}(\psi))=
	\{f \in \mathcal{O}(U) \mid \int_U |f|^2 e^{-\psi} dV < \infty\}
\end{equation}
where $U$ is an open subset of $M$.
This sheaf is also coherent and satisfies the vanishing theorem of Nadel type.
In terms of this formulation, Theorem \ref{thm:Nadel_MIS} can be written as follows.
Let $\{\varphi_t\}$ be the solution $\{\varphi_t\}_{0\le t < t_0}$ of (\ref{eq:Monge-Ampere})  which violates a priori estimates.
\begin{theorem}[\cite{demailly-kollar01}]\label{thm:DK_MIS}
Let $M$ be a Fano manifold of dimension $m$.
Let $G$ be a compact subgroup of $\mbox{Aut}(M)$.
Assume that $M$ does not admit a $G$-invariant K\"ahler-Einstein metric.
Let $\gamma\in (m/(m+1), 1)$.
Then there exists a $G$-invariant sequence $\{\varphi_{t_k}\}_{k=1}^\infty$ such that
\begin{itemize}
\item
	$t_k\to t_0$ as $k\to \infty$,
\item
	there exists a limit $\varphi_\infty=\lim_{k\to\infty}(\varphi_{t_k}-\sup_M\varphi_{t_k})$ in $L^1$-topology, which is an $\omega_g$-psh function, and
\item
	$\mathcal{I}(\gamma\varphi_\infty)$ is a $G^{\mathbb{C}}$-invariant proper multiplier ideal sheaf, i.e, $\mathcal{I}(\gamma\varphi_\infty)$ is neither $0$ nor $\mathcal{O}_M$.
\end{itemize}
\end{theorem}
We call $\mathcal{I}(\gamma \varphi_\infty)$ the KE-MIS of exponent $\gamma$.
In the above, one of the important ingredients is that the upper bound of $\gamma$ is strictly smaller than $1$.
To explain this point, we shall state Nadel's vanishing theorem in terms of Demailly-Koll\'ar's formulation.
Let $L$ be a holomorphic line bundle over $M$ with a singular Hermitian metric $h=h_0e^{-\psi}$, where $h_0$ is a smooth Hermitian metric and $\psi$ is a $L^1_{loc}$-function.
Assume that $\Theta_h(L)=\frac{\sqrt{-1}}{2\pi}\partial\bar{\partial}(-\log h_0+\psi)$ is positive definite in the sense of currents, i.e., $\Theta_h(L)\ge \varepsilon \omega_g$ for some $\varepsilon>0$.
Then, in the same spirit of Nadel's vanishing theorem, we have
\begin{equation}\label{eq:DK_vanishing}
	H^i(M, K_M\otimes L \otimes \mathcal{I}(\psi))=0,
	\,\,\, i>0,
\end{equation}
where $K_M$ is the canonical line bundle.
%%% C³(²–ì)
Now let $\varphi$ be an $\omega_g$-psh function on a Fano manifold $M$ with $[\omega_g]=c_1(M)$.
Substitute $K_M^{-1}$ and $\gamma\varphi$ into $L$ and $\psi$ in (\ref{eq:DK_vanishing}) respectively, and assume $h_0$ associates to $\omega_g$.
Since $\omega_{\varphi}=\omega_g+\frac{\sqrt{-1}}{2\pi}\partial\bar{\partial}\varphi \ge 0$, we have
\[
	\Theta_{h_0e^{-\gamma\varphi}}(L)=\gamma\omega_{\varphi}+(1-\gamma)\omega_g
	\ge (1-\gamma)\omega_g
\]
if $\gamma <1$.
This means that the positivity condition for (\ref{eq:DK_vanishing}) with respect to $h_0e^{-\gamma\varphi}$ holds if $\gamma <1$.
Then (\ref{eq:DK_vanishing}) implies
\[
	H^i(M, \mathcal{I}(\gamma\varphi))=0,
	\,\,\,
	i > 0.
\]
Moreover, if the subscheme cut out by $\mathcal{I}(\gamma\varphi)$ is not empty, then
\[
	H^0(M, \mathcal{I}(\gamma\varphi))=0.
\]
Summing up, we get
\begin{lemma}\label{lem:gamma_vanishing}
If there exists a positive constant $\gamma<1$ and an $\omega_g$-psh function $\varphi$ such that $\mathcal{I}(\gamma\varphi)$ is proper, then $\mathcal{I}(\gamma\varphi)$ satisfies (\ref{eq:vanishing_MISheaf}).
In particular, $\mathcal{I}(\gamma\varphi_\infty)$ for $\gamma \in (m/(m+1),1)$ in Theorem \ref{thm:DK_MIS} satisfies (\ref{eq:vanishing_MISheaf}) (and then (\ref{eq:vanishing_MIScheme})).
\end{lemma}
%%%

On the other hand, the lower bound of $\gamma$ in  Theorem \ref{thm:DK_MIS} describes the strength of the singularity of  $\varphi_\infty$.
It is closely related to a holomorphic invariant introduced by Tian \cite{tian87}.
It is often called the $\alpha$-invariant, which is defined by
\begin{equation}\label{eq:alpha-invariant}
	\alpha_G(M):=
	\sup\{\alpha\in \mathbb{R}
	\mid
	\int_M e^{-\alpha(\psi-\sup_M\psi)}\omega_g^m < C_\alpha
	\,\,\,
	\mbox{for all } G\mbox{-invariant } \omega_g\mbox{-psh } \psi
	\}
\end{equation}
where $G\subset \mbox{Aut}(M)$ is a compact subgroup.
If a multiplier ideal sheaf $\mathcal{I}(\gamma\psi)$ with respect to a $G$-invariant $\omega_g$-psh function $\psi$ of exponent $\gamma$ is proper, where $\sup_M\psi=0$, then $\alpha_G(M)\le \gamma$, because $e^{-\gamma\psi}$ is not integrable over $M$.
Conversely, 
\begin{lemma}\label{lem:alpha_mis}
If $\alpha_G(M) <1$, then there exist a positive constant $\gamma\in (0,1)$ and a $G$-invariant $\omega_g$-psh function $\psi$ with $\sup_M\psi=0$ such that $\mathcal{I}(\gamma\psi)$ is proper.
\end{lemma}
Tian gave a sufficient condition for the existence of K\"ahler-Einstein metrics on Fano manifolds in terms of this invariant.
\begin{theorem}[\cite{tian87}]\label{thm:tian_alpha}
If $\alpha_G(M)>m/(m+1)$, then $M$ admits a $G$-invariant K\"ahler-Einstein metric.
\end{theorem}
Using Theorem \ref{thm:tian_alpha}, Tian and Yau \cite{tian-yau87} proved the existence problem of K\"ahler-Einstein metrics on Fano surfaces, i.e., the Fano surfaces obtained from the blow up of $\mathbb{CP}^2$ at $k$ points where $3 \le k \le 8$ admits a K\"ahler-Einstein metric.
Both of the lower bound of $\alpha_G(M)$ and the non-existence of the proper multiplier ideal sheaves satisfying (\ref{eq:vanishing_MIScheme}) give sufficient condition for the existence of K\"ahler-Einstein metrics on Fano manifolds, and they are related directly to each other.
For example, Lemma \ref{lem:gamma_vanishing} and \ref{lem:alpha_mis}, we have
\begin{lemma}
If $\alpha_G(M) <1$, then a $G^\mathbb{C}$-invariant proper multiplier ideal sheaves satisfying (\ref{eq:vanishing_MIScheme}) exists.
\end{lemma}
Although $\alpha_G(M)$ is difficult to compute in general, it is possible to calculate it when $M$ has a large symmetry such cases as  \cite{song05} for toric varieties and \cite{donaldson0803} for the Mukai-Umemura $3$-folds.
On the other hand, there is a local version of the $\alpha_G(M)$-invariant, which is called the \textbf{complex singularity exponent} \cite{demailly-kollar01}.
Let $K\subset M$ be a compact subset and $\psi$ be a $G$-invariant $\omega_g$-psh function on $M$.
Then the complex singularity exponent $c_K(\psi)$ of $\psi$ with respect to $K$ is defined by
\[
	c_K(\psi)
	=
	\sup\{c\ge 0 \mid e^{-c\psi} \,\, \mbox{is } L^1 \,\, \mbox{on a neighborhood of } K\}.
\]
This constant depends only on the singularity of $\psi$ near $K$.
It is obvious that $c_K(\psi) \ge \alpha_G(M)$.
One of the important properties of $c_K(\psi)$ is the semi-continuity with respect $\psi$. 
Let $\mathcal{P}(M)$ be the set of all locally $L^1$ $\omega_g$-psh functions on $M$ with $L^1$-topology.
Then, we have (cf. Effective version of Main Theorem 0.2 in \cite{demailly-kollar01})
\begin{theorem}[\cite{demailly-kollar01}]\label{thm:semiconti}
Let $K\subset M$ be a compact subset of $M$.
Let $\varphi\in\mathcal{P}(M)$ be given.
If $c<c_K(\varphi)$ and $\psi_j\to\varphi$ in $\mathcal{P}(M)$ as $j\to\infty$, then $e^{-c\psi_j}\to e^{-c\varphi}$ in $L^1$-norm over some neighborhood $U$ of $K$.
\end{theorem}
In particular, if $\{\psi_j\}$ satisfies
\[
	\int_M e^{-\gamma\psi_j}dV \to \infty
\]
where $\gamma \in (\gamma_0, 1)$ and $\psi_j\to\varphi$ in $\mathcal{P}(M)$, then $c_M(\varphi)\le \gamma_0$.
This theorem allows us to substitute Theorem \ref{thm:DK_MIS} for Theorem \ref{thm:Nadel_MIS}.
In fact, if the solution $\varphi_t$ of (\ref{eq:Monge-Ampere}) violates a priori $C^0$-estimate at $t=t_0$, by using a Harnack inequality we can show
\[
	\int_M e^{-\gamma(\varphi_t-\sup\varphi_t)}dV \to \infty
	\,\,\,
	\mbox{as } t\to t_0
\]
for any $\gamma\in (m/(m+1),1)$.
In Theorem \ref{thm:Nadel_MIS}, a subsequence of $\{\varphi_t\}$ induces the KE-MIS, which is proper.
On the other hand, Theorem \ref{thm:semiconti} implies that $e^{-\gamma\varphi_\infty}$ is not integrable over $M$ for any $\gamma\in (m/(m+1),1)$, where $\varphi_{\infty}:=\lim_{i\to\infty}(\varphi_{t_i}-\sup\varphi_{t_i})$.
This means that $\varphi_\infty$ induces the KE-MIS in Theorem \ref{thm:DK_MIS}.

The multiplier ideal sheaves in \cite{demailly-kollar01} and the complex singularity exponent can be defined algebraically as follows (cf. \cite{lazarsfeld_book04-2} and \cite{boyer-galicki08} for instance).
Here we consider a smooth variety $M$ of dimension $m$. 
Let $D=\sum a_iD_i$ be a $\mathbb{Q}$-divisor on $M$.
A \textit{log resolution} of $(M,D)$ is a projective birational map $\mu:M'\to M$ with $M'$ smooth such that the divisor
\[
	\mu^*D+\sum_iE_i
\]
has simple normal crossing support.
Assume $D$ is effective and fix a log resolution $\mu$ of $(M,D)$.
Then the multiplier ideal sheaf $\mathcal{I}(M,D)\subset \mathcal{O}_M$ with respect to $D$ is defined by
\begin{equation}\label{eq:definition_alg_MIS}
	\mu_{*}\mathcal{O}_{M'}(K_{M'}-\mu^*(\lfloor K_M+D  \rfloor)),
\end{equation}
where $\lfloor K_M+D  \rfloor$ means the integral part of $K_M+D$.
Remark that $\mathcal{I}(M,D)$ is independent of the choice of $\mu$.
This (algebraic) ideal sheaf corresponds to the following (analytic) multiplier ideal sheaf defined in \cite{demailly-kollar01}.
Take an open set $U\subset M$ so that for each $D_i$ there is a holomorphic function $g_i$ locally defining $D_i$ in $U$.
Let $\varphi_D:=\sum_i 2a_i \log |g_i|$ which is plurisubharmonic on $U$ and define
\begin{equation}\label{eq:MIS_varphi_D}
	\Gamma(U, \mathcal{I}(\varphi_D)):=
	\{
		f\in \mathcal{O}_M(U) \mid
		\frac{|f|^2}{\prod |g_i|^{2a_i} }
		\in L_{\mbox{loc}}^1
	\}
\end{equation}
as before.
For simplicity, we assume that $D=\sum_ia_iD_i$ has simple normal crossing support.
The holomorphic function $f$ satisfies the $L^2$-integrability condition in (\ref{eq:MIS_varphi_D}) if and only if $f$ can be divided by $\prod g^{m_i}$ where $m_i \ge \lfloor  a_i \rfloor$, i.e.,
$\mathcal{I}(\varphi_D)=\mathcal{O}_M(-\lfloor  D \rfloor)$.
Let $\mu:M'\to M$ be a log resolution of $D$.
Then we have 
\[
	\mathcal{I}(M,D)
	=\mu_*\mathcal{O}_{M'}(K_{M'}-\mu^*(\lfloor K_M+D  \rfloor))
	=\mathcal{O}_M(-\lfloor D  \rfloor)
	=\mathcal{I}(\varphi_{D}).
\]
The second equality in the above was proved in Lemma 9.2.19 \cite{lazarsfeld_book04-2}.
We also have 
\begin{equation}\label{eq:KLT_L1}
	\mathcal{I}(\varphi_D)=\mathcal{O}_M
	\iff
	e^{-\varphi_{D}}\in L^1
	\iff
	(M,D) \,\, \mbox{is KLT}.
\end{equation}
Here we say that a pair $(M,D)$ is KLT if and only if 
\[
	\mbox{ord}_E(K_{M'}-\mu^*(\lfloor  K_M+D \rfloor))>-1
\]
for every exceptional divisor $E$ with respect to a log resolution $\mu:M'\to M$.
In particular $(M,D)$ is KLT  is equivalent to that $\mathcal{I}(M,D)=\mathcal{O}_X$.
The equivalent relation (\ref{eq:KLT_L1}) essentially follows from that if $D_i$ is defined by $\{z_i=0\}$ for a local coordinate $\{z_i\}$ then the $L^1$-integrability of $e^{-\varphi_{D}}$ is equivalent to that $a_i<1$ for all $i$ (Proposition 3.20 \cite{kollar95}).
In particular,
\begin{equation}\label{eq:KLT_condition}
	(M,\gamma D) \,\, \mbox{is KLT}
	\iff
	e^{-\gamma\varphi_{D}} \in L^1.
\end{equation}
Remark that this holds for an $(M,D)$ where $D$ does not necessarily have simple normal crossing support (Proposition 3.20 \cite{kollar95}).
By using the KLT condition (\ref{eq:KLT_condition}), we can rephrase Theorem \ref{thm:DK_MIS}.
Assume that a Fano manifold $M$ does not admit a K\"ahler-Einstein metric.
Let $\varphi_t$ be the solution of (\ref{eq:Monge-Ampere}) where $t\in [0,t_0)$.
As explained before, by taking a subsequence of $\{\varphi_{t_j}-\sup_M\varphi_{t_j}\}$, there exists a limit $\varphi_\infty$ in $L^1$-topology, which is an $\omega_g$-psh function, such that $e^{-\gamma\varphi_\infty} \not\in L^1$ for all $\gamma \in (m/(m+1), 1)$.
Since an approximation theorem in \cite{demailly-kollar01} implies that any $\omega_g$-psh function can be approximated by an $\omega_g$-psh function formed of $\log(\sum_i|f_i|^2)$ where all $f_i$ are holomorphic functions, we can replace the above $\varphi_\infty$ by an $\omega_g$-psh function formed of $\frac{2}{s}\log|\tau_s|$ where $\tau_s\in H^0(M,K^{-s}_M)$ for sufficiently large $s$.
That is to say, there exist a sufficiently large integer $s$ and $\tau_s\in H^0(M,K^{-s}_M)$   such that $e^{-2\gamma\frac{1}{s}\log|\tau_s|}=|\tau_s|^{-\frac{2\gamma}{s}}\not\in L^1$ for all $\gamma\in (m/(m+1), 1)$.
Here $|\cdot|$ is the induced Hermitian metric on $K^{-s}_M$ with respect to the K\"ahler metric $g$.
Hence we have the following theorem.
Remark that the original result holds for orbifolds, but for simplicity we assume that $M$ is smooth in this paper.
\begin{theorem}[Theorem 20 \cite{kollar95}, Theorem 5.2.16 \cite{boyer-galicki08}]
Let $M$ be a Fano manifold of dimension $m$ and $G$ be a compact subgroup of $\mbox{Aut}(M)$.
Assume that there is an $\varepsilon>0$ such that a pair $(M, \frac{m+\varepsilon}{m+1}D)$ is KLT for every $G$-invariant effective divisor $D$ which is numerically equivalent to $K^{-1}_M$.
Then $M$ has a $G$-invariant K\"ahler-Einstein metric.
\end{theorem}
The complex singularity exponent can be also defined algebraically as follows, which is called the \textit{log canonical thresholds} (cf. Appendix in \cite{cheltsov-shramov-demailly0806}).
Let $Z\subset M$ be a closed subvariety.
For an effective $\mathbb{Q}$-Cartier divisor $D$ on $M$, the log canonical threshold of $D$ along $Z$ is defined by
\[
	\mbox{lct}_Z(M,D):=
	\sup\{
		\lambda \in \mathbb{Q} \mid \,\,
		\mbox{the pair } (M, \lambda D) \,
		\mbox{is log canonical along }Z 
	\}.
\]
Here, the pair $(M, D)$ is called log canonical along $Z$ if  $\mathcal{I}(M,(1-\varepsilon)D)$ is trivial in a neighborhood of every point $x\in Z$ for all $0<\varepsilon<1$.
For instance, let us consider  a simple case.
Let $M$ be a Fano manifold and $\sigma \in H^0(M, K^{-l}_M)$.
Let $\psi_\sigma$ be an $\omega$-psh function defined by $\psi_\sigma(z)=\frac{1}{l}\log|\sigma(z)|$.
Let $D_\sigma$ be the associated divisor with $\sigma$ and $Z$ be a closed subvariety in $M$.
In this case, $\mbox{lct}_Z(M, \frac{1}{l}D)$ is the same as $c_Z(\psi_\sigma)$.
The log canonical threshold plays an important role in the studies of the multiplier ideal (sheaves)  in algebraic geometry (cf. \cite{lazarsfeld_book04-2}).
Hence we could expect that the complex singularity exponent with respect to the limit $\varphi_\infty$ in Theorem \ref{thm:DK_MIS} has something to do with the existence of K\"ahler-Einstein metrics although it is not clear at the moment.

To find K\"ahler-Einstein metrics on Fano manifolds, there is another way instead of solving (\ref{eq:Monge-Ampere}), which is the (normalized) K\"ahler-Ricci flow.
The Ricci flow was introduced by R. Hamilton, and on a Fano manifold $M$ with K\"ahler class $c_1(M)$ it is defined by
\begin{equation}\label{eq:KR_flow}
	\frac{d}{dt}\omega_{t}=-\mbox{Ric}(\omega_t)+\omega_t,
	\,\,\,
	\omega_0=\omega_g
\end{equation}
where $t\in [0,\infty)$ and $\omega_t$ is the K\"ahler form of the evolved K\"ahler metric $g_t$.
Remark that (\ref{eq:KR_flow}) is normalized so that the K\"ahler class of $g_t$ is preserved.
The existence and uniqueness of the solution of (\ref{eq:KR_flow}) for $t\in[0,\infty)$ was proved by Cao \cite{cao85}.
If (\ref{eq:KR_flow}) converges in $C^\infty$, the limit is a K\"ahler-Einstein metric.
Then, it is natural to ask whether the results about the multiplier ideal sheaves obtained from the continuity method also hold or not for the K\"ahler-Ricci flow.
The first result of this issue was given by Phong-Sesum-Sturm \cite{pss0611} (see also \cite{phong-sturm0801}).
The equation (\ref{eq:KR_flow}) can be reduced to the equation at the potential level
\begin{equation}\label{eq:KR_flow_potential}
	\frac{d}{dt}\varphi_t=\log(\omega_{t}^m/\omega_g^m)
	+\varphi_t-h_g,
	\,\,\,
	\varphi_0=c
\end{equation}
where $c$ is a constant and $\omega_t=\omega_g+\frac{\sqrt{-1}}{2\pi}\partial\bar{\partial}\varphi_t$.
They gave a necessary and sufficient condition condition for the convergence of $\varphi_t$ as $t \to \infty$.
Their proof consists of the parabolic analogue of Yau's arguments for the elliptic Monge-Amp\`ere equation, the estimates about the K\"ahler-Ricci flow by Perelman 
(cf. \cite{sesum-tian05})
and the result about the Monge-Amp\`ere equations by  Kolodziej (\cite{kolodziej98}, \cite{kolodziej03}).
\begin{theorem}[\cite{pss0611}]\label{thm:PSS_integral}
For a certain appropriate constant $c=c_0$, the convergence of the solution of (\ref{eq:KR_flow_potential}) is equivalent to that there exists $p>1$ such that
\[
	\sup_{t\ge 0} \frac{1}{V}\int_M e^{-p\varphi_t} \omega_g^m
	< \infty.
\]
The convergence is then in $C^\infty$ and exponentially fast.
\end{theorem}
To restate the above theorem in terms of the multiplier ideal sheaves, they introduced 
the sheaf $\mathcal{J}^p$ with respect to a family $\{\psi_t\}_{0\le t <\infty}$ of K\"ahler potentials, which is defined by the presheaf
\begin{equation}\label{eq:MIS_dynamic}
	\Gamma(U, \mathcal{J}^p)=
	\{
		f\in \mathcal{O}(U) \mid
		\sup_{t\ge 0} \int_M |f|^2e^{-p\psi_t} \omega_g^m <\infty
	\}.
\end{equation}
Hence Theorem \ref{thm:PSS_integral} implies
\begin{corollary}[\cite{pss0611}]\label{cor:PSS_MIS}
The K\"ahler-Ricci flow converges if and only if there exists $p>1$ such that $\mathcal{J}^p$ contains the global section $1$.
\end{corollary}
The sheaf $\mathcal{J}^p$ gives a necessary and sufficient condition for the existence of K\"ahler-Einstein metrics and the lower bounds of $p$ is optimal (cf. remarks in \cite{pss0611}), whereas such results are not known for the case of the continuity method.
We emphasize that the sheaf $\mathcal{J}^p$ contains different informations from the multiplier ideal sheaves in Theorem \ref{thm:DK_MIS}, because we do not need the limit of $\{\varphi_t\}$ but the whole of $\{\varphi_t\}$ in order to define $\mathcal{J}^p$.
In fact, in order to get the limit of $\{\varphi_t\}$, we need the appropriate normalization of $\{\varphi_t\}$ as Theorem \ref{thm:DK_MIS}.
In the terminology of the recent paper \cite{siu0902}, the $\mathcal{J}^p$ can be regarded as a dynamic MIS which is similar to the Nadel's formulation rather than a static MIS as the Demailly-Koll\'ar's formulation. 
Now let us consider Theorem \ref{thm:PSS_integral} by using the static MIS instead of the dynamic MIS.
Theorem \ref{thm:PSS_integral} implies that if the normalized K\"ahler-Ricci flow does not converge, then there is a subsequence $\{\varphi_{t_i}\}_i$  of the solution of (\ref{eq:KR_flow_potential}) such that
\[
	\int_M e^{-p(\varphi_{t_i}-\sup\varphi_{t_i})}\omega_g^m \to \infty
\]
as $i\to \infty$ for any $p>1$.
Hence, the limit $\varphi_\infty:=\lim(\varphi_{t_i}-\sup\varphi_{t_i})$ implies the multiplier ideal sheaf $\mathcal{I}(p\varphi_\infty)$, which is proper for any $p>1$.
That is to say, if there is no $G$-invariant $\omega_g$-psh function $\psi$ such that $\mathcal{I}(p\psi)$ is proper for any $p\in (1,+\infty)$, then $M$ admits a $G$-invariant K\"ahler-Einstein metric.
More precisely,
\begin{theorem}\cite{pss0611}\label{thm:PSS_DK_version}
Let $M$ be a Fano manifold.
Let $G\subset \mbox{Aut}(M)$ be a compact subgroup.
Assume that $M$ does not admit K\"ahler-Einstein metrics.
Let $p\in (1,\infty)$ and $\omega_g \in c_1(M)$.
There is a $G$-invariant subsequence of the solutions $\{\varphi_{k_j}\}_{j\ge 1}$ of (\ref{eq:KR_flow_potential}) such that
\begin{itemize}
\item
there exists the limit $\varphi_\infty=\lim_{j\to\infty} (\varphi_{k_j}-\frac{1}{V}\int_M \varphi_{k_j}\omega_g^m)$ in $L^1$-topology, which is an $\omega_g$-psh function, and
\item
$\mathcal{I}(p\varphi_\infty)$ is a $G^{\mathbb{C}}$-invariant proper multiplier ideal sheaf satisfying
\[
	H^i(M, \mathcal{I}(p\varphi_\infty)\otimes K_M^{-\lfloor p \rfloor})
	=0, \,\,\, \mbox{for all }i\ge 1.
\]
\end{itemize}
\end{theorem}
%%%%
Note that Nadel's vanishing formula (\ref{eq:vanishing_MIScheme}) need not hold for the induced MIS $\mathcal{I}(p\varphi_\infty)$, because $p>1$.
However, this result still has an application.
By using a weaker version of Nadel's vanishing theorem and Corollary \ref{cor:PSS_MIS}, Heier \cite{heier0710ii} proved the existence of K\"ahler-Einstein metrics for certain del Pezzo surfaces with large automorphism group.
%%%%

After \cite{pss0611}, Rubinstein \cite{rubinstein0708} gave an analogous result as Theorem \ref{thm:DK_MIS} for the K\"ahler-Ricci flow by using a static MIS as Demailly-Koll\'ar.
His proof is similar to the case of the continuity method, and makes use of the estimates of Perelman and the uniform Sobolev inequality of the K\"ahler-Ricci flow given by Ye \cite{ye0707} and Zhang \cite{zhangQS0708}, which appeared after \cite{pss0611}, in stead of Kolodziej's theorem.
\begin{theorem}[\cite{rubinstein0708}]\label{thm:Rubinstein_MIS}
Let $M$ be a Fano manifold of dimension $m$.
Let $G$ be a compact subgroup of $\mbox{Aut}(M)$.
Let $\gamma\in (m/(m+1),1)$.
Assume that $M$ does not admit a $G$-invariant K\"ahler-Einstein metric.
Then there is an initial constant $c_0=\varphi_0$ and a $G$-invariant subsequence $\{\varphi_{t_i}\}$ of the solution of (\ref{eq:KR_flow_potential}) such that
\begin{itemize}
\item
there exists the limit $\varphi_\infty=\lim_{t_i\to\infty}(\varphi_{t_i}-\frac{1}{V}\int_M\varphi_{t_i}\omega_g^m)$ in $L^1$-topology, which is an $\omega_g$-psh function, and
\item
	$\mathcal{I}(\gamma\varphi_\infty)$ is a $G^{\mathbb{C}}$-invariant proper multiplier ideal sheaf.
\end{itemize}
\end{theorem}
In this paper, we call the above multiplier ideal sheaf $\mathcal{I}(\gamma\varphi_\infty)$ the \textbf{KRF multiplier ideal sheaf} of exponent $\gamma$.
There are some remarks about the above theorem.
First, the KRF multiplier ideal sheaf is independent of the choice of initial constant $c_0$ of (\ref{eq:KR_flow_potential}) due to the normalization of $\varphi_{t_i}$.
In fact, if we choose another constant $c'_0$ instead of the constant $c_0$ in Theorem \ref{thm:Rubinstein_MIS}, which is the same as Theorem \ref{thm:PSS_integral}, the solution of (\ref{eq:KR_flow_potential}) is given by $\varphi'_t=\varphi_t+(c'_0-c_0)e^t$.
In contrast to this, when we consider the convergence of non-normalized K\"ahler potentials $\{\varphi_t\}$ as Theorem \ref{thm:PSS_integral}, we need to pay attention to the choice of the constant $c_0$.
Second, the normalization in Theorem \ref{thm:Rubinstein_MIS} is equivalent to the one in Theorem \ref{thm:DK_MIS}.
In fact, there is a uniform constant $C$ such that $\sup_M \varphi_t -C \le \frac{1}{V}\int_M\varphi_t \omega^m \le \sup_M\varphi_t$.
Third, $\gamma$ is contained in the interval $(m/(m+1),1)$. This means that the subschemes cut out by $\mathcal{I}(\gamma\varphi_\infty)$ satisfies (\ref{eq:vanishing_MIScheme}) and we can make use of the induced geometric properties.
Fourth, the process to prove Theorem \ref{thm:Rubinstein_MIS} 
%%% C³ (²–ì) ‰º'Ì reference 'ðƒRƒƒ"ƒgƒAƒEƒg'µ'Ü'µ'½.
%in \cite{rubinstein0709} 
%%%
is similar to the case of the continuity method, and 
%%% C³ (²–ì) ‰º'Ì reference 'ð rubinstein0709 '©'ç rubinstein0708 'ɕύX'µ'Ü'µ'½.
the proof in \cite{rubinstein0708} 
%%%
implies immediately that if $\alpha_G(M) >\frac{m}{m+1}$ then the K\"ahler-Ricci flow will converge. (This similarity is pointed out in \cite{chen-wang0809} after \cite{rubinstein0708} too.)

Rubinstein \cite{rubinstein0709} also gave the analogous result of Theorem \ref{thm:PSS_DK_version} in terms of the discretization of the K\"ahler-Ricci flow called ``Ricci iteration."
Given a K\"ahler form $\omega\in c_1(M)$ and a real number $\tau>0$, the time $\tau$ Ricci iteration is defined by the sequence $\{\omega_{k\tau}\}_{k\ge 0}$ satisfying
\begin{equation}\label{eq:ricci_iteration}
	\omega_{k\tau}=
	\omega_{(k-1)\tau}+\tau\omega_{k\tau}-\tau \mbox{Ric}(\omega_{k\tau})
	\,\, \mbox{for } k\in \mathbb{N},
\end{equation}
and $\omega_0=\omega$.
When $\tau=1$, (\ref{eq:ricci_iteration}) is the discretization of (\ref{eq:KR_flow}).
Let $\mathcal{H}_\omega$ be the space of K\"ahler potentials with respect to $(\omega, c_1(M))$
\[
	\mathcal{H}_\omega:
	=
	\{
		\psi\in C^\infty_\mathbb{R} \mid
		\omega_\psi=\omega+\frac{\sqrt{-1}}{2\pi}\partial\bar{\partial}\psi >0
	\}.
\]
Let $h_{\omega_\psi}$ be the Ricci potential with respect to $\omega_\psi$ defined as (\ref{eq:ricci_potential}).
Since $[\omega_{k\tau}]=[\omega_{(k-1)\tau}]$, so (\ref{eq:ricci_iteration}) can be written as the system of complex Monge-Amp\`ere equations
\begin{equation}\label{eq:iteration_equation}
	\omega_{\psi_{k\tau}}^m
	=
	\omega^m e^{h_\omega+\frac{1}{\tau}\varphi_{k\tau}-\psi_{k\tau}}
	=
	\omega_{\psi_{(k-1)\tau}}^m
	e^{(\frac{1}{\tau}-1)\varphi_{k\tau}-\frac{1}{\tau}\varphi_{(k-1)\tau}},
\end{equation}
where $k\in \mathbb{N}$, $\omega_{\psi_{k\tau}}=\omega_{k\tau}$ and $\varphi_{k\tau}:=\psi_{k\tau}-\psi_{(k-1)\tau}$.
\begin{theorem}\cite{rubinstein0709}\label{thm:Rubinstein_MIS_iteration}
Let $M$ be a Fano manifold.
Let $G\subset \mbox{Aut}(M)$ be a compact subgroup.
Assume that $M$ does not admit K\"ahler-Einstein metrics.
Let $\tau=1$. 
Let $\gamma\in (1,\infty)$ and $\omega \in c_1(M)$.
There is a $G$-invariant subsequence of the solutions $\{\psi_{k_j}\}_{j\ge 1}$ of (\ref{eq:iteration_equation}) such that
\begin{itemize}
\item
there exists the limit $\varphi_\infty=\lim_{j\to\infty} (\psi_{k_j}-\frac{1}{V}\int_M \psi_{k_j}\omega^m)$ in $L^1$-topology, which is an $\omega$-psh function, and
\item
$\mathcal{I}(\gamma\varphi_\infty)$ is a $G^{\mathbb{C}}$-invariant proper multiplier ideal sheaf satisfying
\[
	H^i(M, \mathcal{I}(\gamma\varphi_\infty)\otimes K_M^{-\lfloor \gamma \rfloor})
	=0, \,\,\, \mbox{for all }i\ge 1.
\]
\end{itemize}
\end{theorem}

%%% ‰Á•M (²–ì) Tosatti 'Æ Szekelyhidi '̘_•¶'Ö'ÌŽQÆ (Phong-Sturm '̘_•¶'à'ljÁ'µ'Ü'µ'½)
Considering Yau's conjecture, it is also natural to ask how stability conditions in the sense of GIT is related to the convergence of the K\"ahler-Ricci flow.
For example, see \cite{phong-sturmJDG06}, \cite{szekelyhidi0803} and \cite{tosatti0810} for references of this issue. 
%%%%

\section{Direct relationships between multiplier ideal sheaves and the obstruction $F$}

It is conjecture by Yau that the existence of canonical K\"ahler metrics such as K\"ahler-Einstein metrics 
and constant scalar curvature metrics for a given K\"ahler class would be equivalent to stability of manifolds in some sense of Geometric Invariant Theory.
This conjecture is formulated by Tian and Donaldson in terms K-polystability as explained in section 1, and is still open.

This conjecture is an analogue of the so-called Hitchin-Kobayashi correspondence, which was proved by Donaldson, and Uhlenbeck and Yau.
The proof of the direction from stability towards the existence of Hermitian-Einstein metrics was proceeded by constructing subsheaves which violate stability (such sheaves are often called \textit{destabilizing subsheaves}) from the bubble of the Yang-Mills heat flow or the continuity method if a vector bundle does not admit a Hermitian-Einstein metric.
Weinkove  \cite{weinkove07} defined a MIS for each sequence of Hermitian metrics on a holomorphic vector bundle and by using it he proved that the bubble of the Yang-Mills heat flow induces a destabilizing subsheaf.
Hence, as the analogy between the Yau-Tian-Donaldson conjecture and the Hitchin-Kobayashi correspondence, we could expect that the MIS obtained from the continuity method or the K\"ahler-Ricci flow  corresponds to a destabilizing subsheaf in some sense for a Fano manifold with anticanonical polarization, but their relation is not clear at this moment.
This issue leads us to study direct relationships between the multiplier ideal sheaves and the obstruction $F$. 
Such a direct relationship was first pointed out by Nadel in \cite{nadel95}.
Extending Nadel's result is the main purpose of this section.

Up to this point we have not defined the character $F$ explicitly, which we do now.
Let $M$ be an $m$-dimensional Fano manifold with K\"ahler class $c_1(M)$ and $g$ be a K\"ahler metric whose K\"ahler form $\omega_g$ represents $c_1(M)$.
We denote the Lie algebra consisting of all holomorphic vector fields on $M$ by $\mathfrak{h}(M)$.
We define the map $F:\mathfrak{h}(M)\to \mathbb{C}$ by
\[
	F(v):=\int_M vh_g \omega_g^m.
\]
In \cite{futaki83.1}, the first author proved that $F$ is independent of the choice of $g$ and that $F$ is a Lie algebra character of $\mathfrak{h}(M)$.
If $M$ is a K\"ahler-Einstein manifold, then $F$ vanishes on $\mathfrak{h}(M)$ because we can take $h_g\equiv 0$.
Thus the vanishing of $F$ is a necessary condition for the existence of K\"ahler-Einstein metrics, but it is known that it is not sufficient.
For example, in \cite{tian97} Tian gave a counterexample which does not admit K\"ahler-Einstein metrics and have no nontrivial holomorphic vector fields.
So it is reasonable to study relationships between the invariant $F$ and the multiplier ideal sheaves.
First of all, we consider the multiplier ideal sheaves obtained from the continuity method in the sense of Nadel.
Assume that $M$ does not admit K\"ahler-Einstein metrics and that $\mathfrak{h}(M)\neq \{0\}$.
For each nontrivial holomorphic vector field $v$, define
\[
	Z^+(v):=\{
	p\in \mbox{Zero}(v) \mid
	\mbox{Re}(\mbox{div}(v))>0
	\},
\]
where $\mbox{Zero}(v)$ is the zero set of $v$ and $\mbox{div}(v)$ is the divergence of $v$ with respect to some K\"ahler metric $g$, i.e., $\mbox{div}(v)=(\mathcal{L}_v(\omega_g^m))/\omega_g^m$,
$\mathcal{L}_v$ being the Lie derivative along $v$.
Remark that $Z^+(v)$ does not depend on the choice of $g$, although $\mbox{div}(v)$ does.
Since $M$ does not admit a K\"ahler-Einstein metric, the closedness of the set of $t's$ for which the solutions $\{\varphi_t\}$ of (\ref{eq:Monge-Ampere}) 
exist does not hold, that is, the solutions cease to exist at some $t_0 \in (0,1]$.
Then the main result in \cite{nadel95} is as follows.
\begin{theorem}[\cite{nadel95}]\label{thm:MIS_Futaki}
Let $M$, $\mathfrak{h}(M)$ and $\{\varphi_t\}$ be as above.
Let $\mathcal{V}$ be the induced KE-MIS obtained from a subsequence of $\{\varphi_{t_i}\}_i$ where $t_i < t_0$ and $t_i \to t_0$.
Then, for any $v\in \mathfrak{h}(M)$ with $F(v)=0$, the support of $\mathcal{V}$ is not contained in $Z^+(v)$.
\end{theorem}
By using the above theorem, Nadel gave another theoretical approach to show that $\mathbb{CP}^1$ does admit a K\"ahler-Einstein metric.
In fact, if we assume that $\mathbb{CP}^1$ did not admit a K\"ahler-Einstein metric, then $\mathcal{V}$ would be zero dimensional and it would be a single reduced point, which follows from (\ref{list:zero_dim}) of the properties of the multiplier ideal subschemes.
We may assume that $\mathcal{V}=\{z=0\}$ in $\mathbb{CP}^1=\mathbb{C}\cup \{\infty\}$.
Let $v=z\frac{d}{d z}\in \mathfrak{h}(M)$, then $v=0$ and the divergence of $v$ is strictly positive at $z=0$. Hence $\mathcal{V}\subset Z^+(v)$, which is a contradiction.
As far as the authors know, other applications of Theorem \ref{thm:MIS_Futaki} except this example had been unknown until \cite{futaki-sano0711}.

We wish to extend this in several ways.We wish first of all to get some more informations about Fano manifolds, 
secondly to show the existence of MIS for K\"ahler-Ricci solitons, and thirdly to study the MIS arising from the non-convergence of K\"ahler-Ricci flow 
and study the relation between MIS and 
%%% C³ (²–ì) "ñ–Ø•s•Ï—Ê'Ì‹L† f 'ð F 'ɏC³'µ'Ü'µ'½.
$F$.
%%%

We study three types of MIS: first of all 
KE-MIS which is due to Nadel, arising from the failure of solving Monge-Amp\`ere equations for K\"ahler-Einstein metrics by continuity method,
secondly KRS-MIS which arises from the failure of solving Monge-Amp\`ere equations for K\"ahler-Ricci solitons by continuity method
and thirdly KRF-MIS which arises from the failure of convergence of K\"ahler-Ricci flow.

Let $M$ be a Fano manifold, 
$G$ a compact subgroup of $\mathrm{Aut}(M)$, and 
$T^r$ the maximal torus of $G$.
For any $G$-invariant K\"ahler metric $g$ with 
$$\omega_g := \frac{\sqrt{-1}}{2\pi} g_{i\barj}dz^i\wedge d\barz^j \in  c_1(M)$$
consider the Hamiltonian $T^r$-action with 
the moment map $\mu_g : M \to \frak t^{r\ast}$. We normalize it by
$$ \int_M u_X e^h \omega^m = 0$$
where $u_X (p) = \langle \mu(p),X\rangle$ and $\mathrm{Ric}_{\omega} - \om = i\p\bp h$. Note that this normalization is
equivalent to requiring $u_X$ to satisfy
$$ \Delta u_X + Xh + u_X = 0,$$
see \cite{futaki88}. 
For $\xi \in \frak t^r$ put
$$ D^{\le 0}(\xi) := \{ y \in \mu(M)\ |\ <y,\xi>\ \le\ 0\}. $$

\begin{theorem} [\cite{futaki-sano0711}]\label{FS1}\ \ 
Suppose $M$ does not admit a K\"ahler-Einstein metric, and let
$V$ be the support of the KE-MIS.
Let $\xi \in \frak t^r \subset \frak h(M)$
satisfy $F(v_\xi) > 0$ 
where
$v_\xi$ is the holomorphic vector field corresponding to $\xi$. Then 
$$ \mu_g (V) \not\subset D^{\le 0}(\xi)$$
for any $G$-invariant K\"ahler metric  $g$ whose K\"ahler form is in $c_1(M)$.
\end{theorem}
\begin{corollary} 
Let $M$ be the one-point blow-up of $\bfC\bfP^2$. 
Then $V$ is the exceptional divisor. 
\end{corollary}
Note that this $V$ destabilizes
slope stability in the sense of 
Ross-Thomas by a result of and Panov and Ross \cite{panov-ross07}.

Here is the outline of the proof of Theorem \ref{FS1}.
Let $h \in C^\infty(M)$ satisfy $\mathrm{Ric}_g - \omega_g = i\partial\barpartial h$.
Suppose
$$ \frac{\det(g_{i\barj} + \varphi_{i\barj})}{\det(g_{i\barj})} = e^{-t\varphi + h}$$
has solutions only for $t \in [0,t_0),\ t_0 < 1$.
Then we have  an MIS with a support $V$. The following fact is due to Nadel based on earlier estimates by Siu and Tian.

\begin{fact} \label{fact1}
Let $K \subset M - V$ be a compact subset of $M - V$. Then
$$ \int_K \omega_{g_t}^m \to 0$$
as $t \to t_0$.
\end{fact}

\begin{fact}\label{fact2}
$$\mu_g(p) \in D^{\le 0}(\xi) \Longleftrightarrow (\mathrm{div}(v_\xi))(p) \ge 0$$
where
$$\mathrm{div}(v_\xi)(e^h\omega^m) = \mathcal L_{v_\xi}(e^h\omega^m).$$
\end{fact}

\begin{fact}\label{fact3}
$$\frac t{t-1} F(v_\xi) = \int_M \mathrm{div}(v_\xi) \omega_t^m.$$
\end{fact}

By Fact \ref{fact3} and our assumption $F(v_\xi) > 0$, we have for $t \in (\delta, t_0)$ with $t_0 < 1$
$$ \int_M \mathrm{div}(v_\xi)\ \omega_t^m = \frac t{t-1} F(v_\xi) < - C$$
with $C > 0$ independent of $t$.

We seek a contradiction by assuming $\mu_g(V) \subset D^{\le 0}(\xi) = \{\mathrm{div}(v_\xi) \ge 0$\}.
Choose $\epsilon > 0$ small and put
$$ W_\epsilon := \{p\in M | \mathrm{div}(v_\xi)(p) \le - \epsilon\}.$$
Then $W_\epsilon \subset M - V$ compact. Apply Fact \ref{fact1} to $W_\epsilon$ to get
$$ \int_{W_\epsilon} \omega_{g_t}^m \to 0 $$
as $t \to t_0$.

But then
\begin{eqnarray*}
-C \ge \int_M \mathrm{div}(v_\xi) \omega_t^m &=& \int_{M - W_{\epsilon}} \mathrm{div}(v_\xi) \omega_t^m +
\int_{W_{\epsilon}} \mathrm{div}(v_{\xi}) \omega_t^m\\
&\ge& -2\epsilon \mathrm{vol}(M, g)
\end{eqnarray*}
as $t \to t_0$, a contradiction ! This completes the outline of the proof of Theorem \ref{FS1}.

Next, we turn to KRS-MIS. Let $M$ be again a Fano manifold of dimension $m$. Let 
$\omega_g \in c_1(M)$ be  a K\"ahler form and 
$v \in \frak h_r(M)$ a holomorphic vector field in the reductive part $\frak h_r(M)$ of $\frak h(M)$.

\begin{definition} The pair $(g,v)$ is said to be a K\"ahler-Ricci soliton if 
$$\mathrm{Ric}(\omega_g) - \omega_g = \mathcal L_v(\omega_g).$$
(Hence $\Im(v)$ is necessarily a Killing vector field.)
\end{definition}

Start with an initial metric $g^0$ with $\omega_0 := \omega_{g^0} \in c_1(M)$.
$$\mathrm{Ric}(\omega_0) - \omega_0 = i\partial\barpartial h_0,\quad \int_M e^{h_0}\omega_0^m = \int_M \omega_0^m.$$
$$ i_v \omega_0 = i\barpartial \theta_{v,0}, \quad \int_M e^{\theta_{v,0}}\omega_0^m = \int_M \omega_0^m.$$
Consider for $t \in [0,1]$ 
$$ \det (g^0_{i\barj} + \varphi_{ti\barj}) = \det (g^0_{i\barj}) e^{h_0 - \theta_{v,0} - v\varphi_t - t\varphi_t}.$$
The solution for $t=1$ gives the K\"ahler-Ricci soliton.
Zhu \cite{Zhu} has shown that $t=0$ always has a solution. The 
implicit function theorem shows for some $\epsilon > 0$, 
all $t \in [0,\epsilon)$ have a solution. 

Suppose we only have solutions on $[0, t_{\infty}),\ t_{\infty} < 1$.

Let $\theta_{v,g}$ satisfy 
$$ i_v \omega_g = i\barpartial \theta_{v,g}, \quad \int_M e^{\theta_{v,g}}\omega_g^m = \int_M \omega_g^m.$$

\begin{definition} Define $F_v : \frak h(M) \to \bfC$ by 
$$F_v(w) = \int_M w(h_g - \theta_{v,g}) e^{\theta_{v,g}} \omega_g^m.$$
\end{definition}
\noindent
Tian and Zhu \cite{tian-zhu} showed that this $F_v$ is independent of $g$ with $\omega_g \in c_1(M)$. 

\begin{theorem}[Tian-Zhu \cite{tian-zhu}]\label{v}   There exists a unique $v \in \frak h_r(M)$ such that 
$$ F_v(w) = 0 \mathrm{\ for\ all\ }w \in \frak h_r(M).$$
\end{theorem}
\noindent
We take $v$ to be the one chosen in the Theorem \ref{v}.

\begin{theorem}[\cite{futaki-sano0711}] \label{FS2}Let $K$ be the compact subgroup such that $\frak k \otimes \bfC = \frak h_r (M)$.
Let $v$ be the one chosen in the Theorem \ref{v}.
Suppose there is no KRS. Then we get an MIS and its support $V_s$ satisfies
$$ V_v \not\subset Z^+(\mathrm{grad}'w) \mathrm{\ for\ } \forall v \in \frak h_r(M).$$
\end{theorem}

Just as Nadel applied Theorem \ref{thm:MIS_Futaki} to prove the existence of K\"ahler-Einstein metric on $\bfC\bfP^1$, 
we can apply Theorem \ref{FS2} to prove the existence of KRS on the one point blow-up of $\bfC\bfP^2$.

Next we consider KRF-MIS. As mentioned in section 2, there are two approaches to KRF-MIS, one by Phong, Sesum and Sturm \cite{pss0611},
and the other by Rubinstein \cite{rubinstein0708}. Here, we consider the one considered by Rubinstein. So, 
one gets an MIS from the failure of convergence of normalized K\"ahler-Ricci flow:
\begin{equation}\label{eq:KR_flow_again}
	\frac{\partial g}{\partial t} = - \mathrm{Ric}(g) + g.
\end{equation}
If we put $g_{ti\barj} = g_{i\barj} + \varphi_{ti\barj}$ the Ricci flow is equivalent to
$$ \frac{\partial \varphi_t}{\partial t} = \log \frac{\det(g_{i\barj} + \varphi_{ti\barj})}{\det(g_{i\barj})} 
+ \varphi_t - h_0$$
$$ \varphi_0 = c_0$$
Rubinstein modified Phong-Sesum-Sturm's MIS using the idea of Demailly-Koll\'ar:
$$ \varphi_t - \int_M \varphi_t\,\omega^m \longrightarrow \varphi_{\infty} \quad\mathrm{almost\ psh } $$
as $t \to \infty$.
Let $V_\gamma$ be the MIS for $\psi = \gamma \varphi_{\infty}$, $\gamma \in (\frac m{m+1}, 1)$, 
defined by
$$ \Gamma (U, \mathcal I(\psi)) = \{\ f\ \in \mathcal O_M(U)\ |\ \int_U\ |f|^2\ e^{-\psi}\ \omega_g^m < \infty\}.$$
This MIS satisfies
$$ H^q(M, \mathcal I(\psi)) = 0 \mathrm{\quad for\quad } \forall q > 0.$$
In general, it seems to be difficult to calculate (the support) of KRF-MIS.
However, under some assumptions, it becomes computable.
To explain it, we recall the following two results.
Let $M$ be a toric Fano manifold of dimension $m$, on which the algebraic torus $T_\bfC = (\bfC^\ast)^m$ acts. 
Let $T_\bfR = T^m$ be the real torus and $\mathfrak t_\bfR$ its Lie algebra. We further put
$N_\bfR = J\mathfrak t_\bfR$. Let 
$W(M) = N(T_\bfC)/T_\bfC$ be the Weyl group. 

\begin{theorem}[Wang-Zhu \cite{wang-zhu}]\label{thm:wang-zhu} There exists a K\"ahler-Ricci Soliton $(g_{\mathrm{KRS}}, v_{\mathrm{KRS}})$.
\end{theorem}

Here we assume that $K$ denotes a maximal compact subgroup of the reductive part of $\mathrm{Aut}(M)$ and $K_{v_{\mathrm{KRS}}}$ denotes the one-parameter subgroup of $K$ generated by the imaginary part of $v_{\mathrm{KRS}}$.
Then, 
\begin{theorem}[Tian-Zhu \cite{tian-zhu07}]\label{thm:tian-zhu}
Let $M$ be a (not necessarily toric) Fano manifold which admits a K\"ahler-Ricci soliton $(g_{\mathrm{KRS}}, v_{\mathrm{KRS}})$.
Then, any solution $g_t$ of (\ref{eq:KR_flow_again}) will converge to $g_{\mathrm{KRS}}$ in the sense of Cheeger-Gromov if the initial K\"ahler metric is $K_{v_{\mathrm{KRS}}}$-invariant.
\end{theorem}
Combining Theorem \ref{thm:wang-zhu} and Theorem \ref{thm:tian-zhu}, we find that the flow (\ref{eq:KR_flow_again}) always converges to a K\"ahler-Ricci soliton in the sense of Cheeger-Gromov on toric Fano manifolds.
This fact suggests us a possibility to understand the asymptotic behavior of $g_t$ along (\ref{eq:KR_flow_again}) and to get some information about KRF-MIS from data of K\"ahler-Ricci solitons.  
In fact, the second author proved
\begin{theorem}[\cite{sano0811}] 
Suppose that the fixed point set $N_\bfR^{W(M)}$ of  the Weyl group $W(M)$ on $N_\bfR$ is one dimensional. Let 
$\sigma_t = \exp(tv_{\mathrm{KRS}})$ be the one parameter group of transformations generated by $v_{\mathrm{KRS}}$, $0 < \gamma <1$ and 
$\omega$ a $T_\bfR$-invariant K\"ahler form in $c_1(M)$.
Then the support of Rubinstein's KRF-MIS of exponent $\gamma$ 
is equal to the support of 
the MIS of exponent $\gamma$  
obtained from the K\"ahler potentials of $\{(\sigma_t^{-1})^\ast\omega\}$.
\end{theorem}
Remark that the assumption of $N_\bfR^{W(M)}$ is constrained and it would be expected to be removed. 
Using the above theorem, the second author computed the support of KRF-MIS for various $\gamma$ on some examples.
For example, we can prove
\begin{corollary}
Let $M$ be the blow up of $\mathbb{CP}^2$ at $p_1$ and $p_2$.
Let $E_1$ and $E_2$ be the exceptional divisors of the blow up, and $E_0$ be the proper transform of $\overline{p_1p_2}$ of the line passing through $p_1$ and $p_2$.
Then, the support of KRF-MIS on $M$ of
exponent $\gamma$ is 
\[
	\left\{\begin{array}{cc}
		\cup_{i=0}^{2} E_i & \mbox{for } \gamma \in (\frac{1}{2}, 1),  \\
		E_0 & \mbox{for } \gamma \in (\frac{1}{3}, \frac{1}{2}).
	\end{array}\right.
\]
\end{corollary}
It would be interesting to consider a relationship between destabilizing test configurations and the pair of the support of KRF-MIS and its exponent.
%%%


\begin{thebibliography}{99}

\bibitem{aubin76}T.~Aubin : Equations du type de Monge-Amp\`ere sur les 
vari\'et\'es k\"ahl\'eriennes compactes, C. R. Acad. Sci. Paris, {\bf 283}, 
119-121 (1976)



\bibitem{boyer-galicki08} C.P. Boyer and K. Galicki : Sasakian geometry, (Oxford Mathematical Monographs., 2008).

\bibitem{calabi85}E.~Calabi : Extremal K\"ahler metrics II, Differential
geometry and complex analysis, (I. Chavel and H.M. Farkas eds.),  
95-114, Springer-Verlag, Berline-Heidelberg-New York,
(1985).

\bibitem{cao85}
H.D. Cao :  Deformation of K\"ahler metrics to K\"ahler-Einstein metrics on compact K\"ahler manifolds, Invent. Math. 81 (1985) 359-372. 


\bibitem{cheltsov-shramov-demailly0806} I. Cheltsov and C. Shramov : Log canonical thresholds of smooth Fano threefolds (with an appendix by J.P. Demailly), arXiv:0806.2107 (2008).

\bibitem{chentian04}X.X.~Chen and G.~Tian : Geometry of K\"ahler metrics and 
foliations by holomorphic discs, Publ. Math. Inst. Hautes \'Etudes Sci. No. 107 (2008), 1-107.
math.DG/0409433

\bibitem{chen-wang0809} X.X. Chen and B. Wang : Remarks on K\"ahler Ricci flow, arXiv:0809.3963 (2008).

\bibitem{CFO} K.~Cho, A.~Futaki and H.~Ono : 
 Uniqueness and examples of toric Sasaki-Einstein manifolds, 
 Comm. Math. Phys., 277 (2008), 439-458, math.DG/0701122


\bibitem{demailly-kollar01} J.P. Demailly and J. Koll\'ar :
Semi-continuity of complex singularity exponents and K{\"a}hler-{Einstein} metrics on Fano orbifolds, Ann. Sci. {\'E}cole Norm. Sup. (4) 34, no.4 (2001) 525--556.

\bibitem{donaldson02}S.K.~Donaldson : Scalar curvature and stability of toric
varieties, J. Differential Geometry, 62(2002), 289-349.

\bibitem{donaldson05}S.K.~Donaldson : Lower bounds on the Calabi functional,
J. Differential Geometry, 70(2005), 453-472.

\bibitem{donaldson0803} S.K. Donaldson : K\"ahler geometry on toric manifolds, and some other manifolds with large symmetry, 
Handbook of geometric analysis. No. 1, 29--75, Adv. Lect. Math. (ALM), 7, Int. Press, Somerville, MA, 2008. arXiv:0803.0985 (2008).

\bibitem{donaldson97}S.K.~Donaldson : Remarks on gauge theory, complex geometry 
and four-manifold topology, in 'Fields Medallists Lectures' (Atiyah, Iagolnitzer eds.), 
World Scientific, 1997, 384-403.

\bibitem{donkro}S.K.~Donaldson and P.B.~Kronheimer : The geometry of four manifolds, Oxford Mathematical Monographs, Claren Press, Oxford, 1990.

\bibitem{fujiki92}A.~Fujiki : 
Moduli space of polarized algebraic manifolds and K\"ahler metrics, Sugaku Expositions,
5(1992), 173-191.

\bibitem{futaki83.1}A.~Futaki : 
An obstruction to the existence of Einstein K\"ahler metrics, Invent. 
Math. 73, 437-443 (1983).

\bibitem{futaki83.2}A.~Futaki : On compact K\"ahler manifolds of constant scalar curvature, 
Proc. Japan Acad., Ser. A, 59, 401-402 (1983).

\bibitem{futaki88}A.~Futaki : K\"ahler-Einstein metrics and integral invariants,
Lecture Notes in Math., vol.1314, Springer-Verlag, Berline-Heidelberg-New York,(1988).

\bibitem{futaki05}A.~Futaki : Stability, integral invariants and canonical K\"ahler metrics, 
Proc. 9-th Internat. Conf. on Differential Geometry and its Applications, 2004 Prague, (eds. J. Bures et al),
2005, 45-58, MATFYZPRESS, Prague.

\bibitem{FOW}A.~Futaki,  H.~Ono and G.~Wang : 
 Transverse K\"ahler geometry of Sasaki manifolds and toric Sasaki-Einstein manifolds,
to appear in J. Differential Geom., math.DG/0607586.


\bibitem{futaki-sano0711} A. Futaki and Y. Sano : Multiplier ideal sheaves and integral invariants on toric Fano manifolds, arXiv:0711.0614 (2007).


\bibitem{heier0710i} G. Heier : Convergence of the K\"ahler-Ricci flow and multiplier ideal sheaves on Del Pezzo surfaces, 
to appear in Michigan Math. J., arXiv:0710.5725 (2007).


\bibitem{heier0710ii} G. Heier : Existence of K\"ahler-Einstein metrics and multiplier ideal sheaves on Del Pezzo surfaces, 
to appear in Math. Zeit., arXiv:0710.5724.

\bibitem{kollar95} J. Koll\'ar : Singularities of pairs. (Algebraic Geometry---Santa Cruz 1995), Proc. Sympos. Pure Math., 62, part 1 (1995) 221--287.


\bibitem{kolodziej98} S. Kolodziej : The complex Monge-Amp\`ere equation, Acta Math. 180, no.1 (1998) 69--117.

\bibitem{kolodziej03} S. Kolodziej : The Monge-Amp\`ere equation on compact K\"ahler manifolds, Indiana Univ. Math. J. 52, no.3 (2003) 667--686.

\bibitem{lazarsfeld_book04-2} R. Lazarsfeld : Positivity in Algebraic Geometry II: Positivity for vector bundles, and multiplier ideals, (Springer-Verlag, 2004).

\bibitem{mabuchi0812}T.~Mabuchi : 
K-stability of  constant scalar curvature polarization, 
arXiv:0812.4093

\bibitem{matsushima57}Y.~Matsushima : Sur la structure du groupe 
d'hom\'eomorphismes d'une certaine vari\'et\'e kaehl\'erienne, Nagoya
Math. J., {\bf 11}, 145-150 (1957).


\bibitem{nadel90}
A.M. Nadel : Multiplier ideal sheaves and K\"ahler-Einstein metrics of positive scalar curvature, Ann. of Math. (2) 132, no.3 (1990) 549--596.

\bibitem{nadel95}
A.M. Nadel : Multiplier ideal sheaves and Futakifs invariant, Geometric Theory of Singular Phenomena in Partial Differential Equations (Cortona, 1995), Sympos. Math., XXXVIII, Cambridge Univ. Press, Cambridge, 1998. (1995) 7--16.

\bibitem{panov-ross07}D.~Panov and J. Ross : Slope Stability and Exceptional Divisors of High Genus, Math. Ann. 343 (2009), no. 1, 79--101.
arXiv:0710.4078v1 [math.AG].


\bibitem{paultianCM1}S.T.~Paul and G.~Tian : CM Stability and the Generalized Futaki Invariant I., math.AG/0605278, 
2006. 

\bibitem{paultianCM2}S.T.~Paul and G.~Tian : CM Stability and the Generalized Futaki Invariant II. (To appear in Asterisque),
 math.AG/0606505, 2006. 


\bibitem{phong-sturmJDG06} D.H. Phong and J. Sturm : On stability and the convergence of the K\"ahler-Ricci flow, J. Differential Geom.  72  (2006),  no. 1, 149--168.

\bibitem{phong-sturm0801} D.H. Phong and J. Sturm : Lectures on stability and constant scalar curvature, arXiv:0801.4179 (2008).


\bibitem{pss0611} D.H. Phong, N. Sesum and J. Sturm :  Multiplier ideal sheaves and the K\"ahler-Ricci flow, Comm. Anal. Geom. 15, no. 3 (2007), 613--632.

\bibitem{rossthomas06}J.~Ross and R.P.~Thomas : An obstruction to the existence of constant scalar curvature K\"ahler metrics,
J. Differential Geom. 72 (2006), no. 3, 429--466.



\bibitem{Rubin0706}Y.A.~Rubinstein : The Ricci iteration and its applications, C. R. Acad. Sci. Paris, Ser. I, 345 (2007), 445-448, arXiv:0706.2777.

\bibitem{rubinstein0708}Y.A.~Rubinstein : On the construction of Nadel multiplier ideal sheaves and the limiting behavior of the Ricci flow, to appear in Transact. Amer. Math. Soc., arXive:math/0708.1590.

\bibitem{rubinstein0709}Y.A.~Rubinstein : Some discretizations of geometric evolution equations 
and the Ricci iteration on the space of K\"ahler metrics, Adv. Math. 218 (2008), no. 5, 1526--1565.  arXive:math/0709.0990.

\bibitem{sano0811}Y.~Sano : Multiplier ideal sheaves and the K\"ahler-Ricci flow 
on toric Fano manifolds with large symmetry, preprint, arXiv:0811.1455.

\bibitem{sesum-tian05} N. Sesum and G. Tian : Bounding scalar curvature and diameter along the K\"ahler-Ricci flow (after Perelman), 
J. Inst. Math. Jussieu 7 (2008), no. 3, 575--587.

\bibitem{siu88} Y-T. Siu : The existence of K\"ahler-Einstein metrics on manifolds with positive anticanonical line bundle and a suitable finite symmetry group, Annals of Math. 127 (1988) 585--627.

\bibitem{siu0902} Y-T. Siu : Dynamical multiplier ideal sheaves and the construction of rational curves in Fano manifolds, arXiv:0902.2809 (2009).


\bibitem{song05} 
J. Song : The $\alpha$-invariant on toric Fano manifolds, Amer. J. Math. 127, No.6 (2005) 1247--1259.

\bibitem{stoppa0803}J.~Stoppa : 
K-stability of constant scalar curvature K\"ahler 
manifolds, arXiv:0803.4095

\bibitem{szekelyhidi0803} G.~Sz\'ekelyhidi : The K\"ahler-Ricci flow and K-polystability, arXiv:0803.1613 (2008).


\bibitem{tian87}
G. Tian : On K\"ahler-Einstein metrics on certain K\"ahler manifolds with $C_1(M)>0$, Invent. Math. 89, no.2 (1987) 225--246.


\bibitem{tian97} G. Tian, K{\"a}hler-{Einstein} metrics with positive scalar curvature, Invent. Math. 130, no.1 (1997) 1--37.



\bibitem{tian-yau87} G. Tian and S.T. Yau,K\"ahler-Einstein metrics on complex surfaces with $C_1(M)$ positive, Comm. Math. Phys. 112 (1987) 175--203.



\bibitem{tian-zhu}
G. Tian and X. Zhu : A new holomorphic invariant and uniqueness of K\"ahler-Ricci solitons, 
Comment. Math. Helv 77, No.2 (2002) 297--325.

\bibitem{tian-zhu07}
G. Tian and X. Zhu : Convergence of K\"ahler-RIcci flow, J. Amer. Math. Soc., 20, No.3 (2007), 675--699.

\bibitem{tosatti0810} V.~Tosatti : K\"ahler-Ricci flow on stable Fano manifolds, arXiv:0810.1895 (2008).

\bibitem{xwang04}X.-W.~Wang : 
Moment maps, Futaki invariant and stability of projective manifolds, 
Comm. Anal. Geom. 12 (2004), no. 5, 1009--1037.

\bibitem{wang-zhu}
X.J. Wang and X. Zhu : K\"ahler-Ricci solitons on toric manifolds with positive first Chern class, Adv. Math. 188, No.1 (2004) 87--103.


\bibitem{weinkove07} B. Weinkove, A complex Frobenius theorem, multiplier ideal sheaves and Hermitian-Einstein metrics on stable bundles, Trans. Amer. Math. Soc. 359, no.4 (2007) 1577--1592.


\bibitem{yau78}S.-T.Yau : On the Ricci curvature of a compact K\"ahler
manifold and the complex Monge-Amp\`ere equation I, Comm. Pure Appl.
Math. 31(1978), 339-441.

\bibitem{yau93}S.-T.~Yau : Open problems in Geometry, Proc. Symp. Pure Math. 54 (1993) 1-28.

\bibitem{ye0707} R. Ye : The logarithmic Sobolev inequality along the Ricci flow, arXiv:0707.2424 (2007).

\bibitem{zhangQS0708} O.S. Zhang : A uniform Sobolev inequality under Ricci flow, arXiv:0706.1594 (2007).

\bibitem{Zhu}X.~Zhu, K\"ahler-Ricci soliton type equations on compact complex 
manifolds with 
$C_1(M) > 0$, J. Geom. Anal., 10 (2000), 759-774.


\end{thebibliography}
\end{document}